\documentclass{article}
\usepackage{amsmath}
\usepackage{chicago}
\usepackage{epsfig}

\newtheorem{theorem}{Theorem}

\newtheorem{corollary}[theorem]{Corollary}

\input{tcilatex}

\begin{document}

\title{Lattice Option Pricing By Multidimensional Interpolation}
\author{Vladislav Kargin\thanks{%
Cornerstone Research, 599 Lexington Avenue floor 43, New York, NY 10022;
skarguine@cornerstone.com}}
\maketitle

\begin{abstract}
This note proposes a method for pricing high-dimensional American options
based on modern methods of multidimensional interpolation. The method allows
using sparse grids and thus mitigates the curse of dimensionality. A
framework of the pricing algorithm and the corresponding interpolation
methods are discussed, and a theorem is demonstrated that suggests that the
pricing method is less vulnerable to the curse of dimensionality. The method
is illustrated by an application to rainbow options and compared to Least
Squares Monte Carlo and other benchmarks.
\end{abstract}

\begin{quotation}
The fundamental problem of options theory is the valuation of hybrid,
non-linear securities, and options theory is an ingenious but glorified
method of interpolation.

Emanuel Derman ``A guide for the perplexed quant''
\end{quotation}

\section{Introduction}

Lattice option pricing\footnote{%
Invented by \citeN{cox_ross_rubinstein79}.} is very popular among
practitioners because in one-dimensional situations it is straightforward to
implement and has a transparent interpretation. In multiple dimensions it
has two serious drawbacks: the need to build recombining trees and the curse
of dimensionality. If branches of the tree do not recombine then the number
of nodes grows exponentially with the number of time steps. Similarly, the
number of points in a regular grid grows exponentially with the space
dimension. For example, to approximate every point in the nine-dimensional
hypercube with 10\% precision we need one billion points. With twenty
factors, the argument runs, experiments with IBM's fastest supercomputer
will quickly convince us that the lattice method of pricing is impractical
with many dimensions.\footnote{%
The speed of this supercomputer is around $10^{13}$ operations per second,
so an optimistic estimate of the execution time is $10^{7}$ seconds or about 
$100$ days. For pricing algorithms that generalize lattice pricing to
multiple dimensions see \citeN{boyle88}, \citeN{boyle_evnine_gibbs89}, %
\citeN{madan_milne_shefrin89}, \citeN{kamrad_ritchken91}, and %
\citeN{mccarthy_webber01}.}

This argument is usually well taken but essentially wrong. First, trees may
very well have non-recombining branches and a moderate number of nodes. We
only need to leave some nodes without descendants and to interpolate values
on these nodes. Second, in the case of many dimensions, we can use irregular
grids with a moderate number of points provided that we can interpolate the
value function faithfully from its values on the irregular grid.

A suspicion may arise that this just shifts the computational burden to the
interpolation problem and that the curse of dimensionality will remain as
dangerous as it was before. Indeed, traditional approximation theory says
that accurate approximation of a general function needs a number of grid
points that is exponential in the space dimension.

Does this invalidate the idea of pricing by interpolation? No, because we
deal with specific classes of functions which may be approximated better
than an arbitrarily chosen function. As the simplest example, consider
linear functions. They can be recovered from values on only $d+1$ points in $%
d$-dimensional space. Recall also the classical Kotelnikov-Shannon sampling
theorem that says that a function with limited bandwidth of its Fourier
transform can be completely recovered from the values it takes on a discrete
sample of points. These examples suggest that for certain classes of
functions the problem of approximation may be satisfactorily solved even in
multiple dimensions. Indeed, this hope is being realized in recent and
continuing work on non-linear multidimensional approximation. This paper
aims to apply these new ideas to option pricing.

The essential idea is to choose the approximating functions adaptively.
First, the researcher chooses a large class of functions that are easy to
compute and that are suitable to the problem at hand. Second, the researcher
allows the data to select the functions that are most fitting as a basis of
approximation. As a result, the approximation is well adapted to the
properties of the given function. The main practical and theoretical tasks
are automating this procedure and exploring its convergence properties.

The idea of using approximations for option pricing is not entirely new.
Recently, it was implemented by \citeN{longstaff_schwartz01} in their
``simple but powerful'' Least Squares Monte Carlo algorithm.\footnote{%
See also related work by \citeN{tilley93}, \citeN{barraquand_martineau95}, %
\citeN{carriere96}, \citeN{broadie_glasserman97a}, %
\citeN{broadie_glasserman97b}, and \citeN{raymar_zwecher97}.} Their method
is based on the Monte Carlo pricing method and they assume that the
researcher can guess good basis functions for approximations. %
\citeN{tsitsiklis_roy01} describe a similar method and analyze its
convergence properties. The cardinal distinction of the method described in
this paper is that it suggests the adaptive choice of the approximating
basis. This adaptive way of approximation allows construction of an
algorithm which is universally applicable to a wide range of possible
options.

As an additional benefit, the approximation provided by the new algorithm
allows easy pricing everywhere in factor space. In contrast, both the
standard lattice and Monte Carlo methods produce option values for only one
combination of factors. This benefit is especially useful if there is a need
to visualize the dependence of option price on factors, or to compute hedge
factors - sensitivities of the option value with respect to changes in
factors.

The remainder of the paper is organized as follows. Section~2 briefly
reviews the problem of option pricing and relates it to the problem of
approximation. Section~3 outlines the framework of the algorithm,
illustrates it with a simple example, and briefly describes why adaptive
approximation is likely to break the curse of dimensionality. Section~4
explaines how the approximation can be used to find lower and upper bounds
on the option value. Section 5 applies the method to a set of benchmark
options and compare the results with the results in the literature.
Section~6 concludes.

\section{Option Pricing and Interpolation}

The problem is as follows. Let the price of a derivative security depend on $%
N$ factors that follow a specified diffusion process. The derivative is of
the American type and so can be exercised at any time. Assume also that the
price of the derivative is not path-dependent. Then by the dynamic
replication argument, the value of the derivative satisfies the familiar%
\footnote{%
see, for example, \citeN{hull99} or \citeN{wilmott_howson95}.} partial
differential equation :%
\begin{equation}
\frac{\partial f}{\partial t}+r\sum_{i}x_{i}\frac{\partial f}{\partial x_{i}}%
+\frac{1}{2}\sum_{i,k}v_{ik}x_{i}x_{k}\frac{\partial ^{2}f}{\partial
x_{i}\partial x_{k}}=rf.  \label{pde}
\end{equation}%
Here $f(x,t)$ is the value of the derivative at time $t$ if the vector of
factors is $x.$

One way to solve this equation is to write it in an integral form using a
certain measure $\mu $ over the space of Brownian motion paths. This is the
Feynman-Kac representation\footnote{%
see \citeN{kac49} or \citeN{karatzas_shreve91}.} of PDE (\ref{pde}) solution
as an integral from a functional of Brownian motion paths:%
\begin{equation}
f(x,0)=\sup_{\tau }\int_{x(t)\in \mathcal{W}}\pi (x(\tau ),\tau )e^{-rt}d\mu
(x(t)).
\end{equation}%
Here $\tau $ is a stopping time, $\pi (x(\tau ),\tau )$ is the payoff at
time $\tau $ if the factors are $x(\tau ),$ and $\mathcal{W}$ is the space
of paths of the Wiener process.

We can write the Wiener measure $\mu $ as a limit over a sequence of time
discretizations with Gaussian transition probabilities. Then, each of the
discrete time problems can be solved recursively through the Bellman
equation that relates the current option value to the values at the next
time stage:%
\begin{equation}
f(x,t)=\max \left\{ \pi (x(t),t);\int_{y\in X}f(y,t+\Delta t)e^{-r\Delta t}d%
\widetilde{\mu }(x,y,t)\right\} .  \label{bellman}
\end{equation}%
Here $\Delta t$ is a time interval, $d\widetilde{\mu }(x,y,t)$ is the
probability of transition from $x$ to $\ y$ consistent with the Wiener
measure $d\mu ,$ and $X$ is the space of factors.

The next step - crucial for our analysis - is to discretize the equation
over space. \ This means choosing a grid $G\subset X$ and a suitable
approximation for $d\widetilde{\mu }$. Here is where difficulties begin. The
standard tree methods use a regular -- usually cubic -- lattice, and specify
probabilities of transitions from each lattice point to nearby lattice
points so as to match the covariance matrix of the continuous process.

In multidimensional situations the number of points in regular lattices is
prohibitively large. Therefore we have to use an irregular grid with large
gaps. Transitions over such a grid are unlikely to approach the Wiener
process uniformly.\footnote{%
See, however, \citeN{berridge_schumacher04} for encouraging advances in this
direction.} So what to do? One solution is to use spatial interpolation as
in the following formula:

\begin{equation}
\widehat{f}(x,t)=\max \left\{ \pi (x(t),t),\int_{y\in \mathcal{D}(x)}%
\widehat{f}(y,t+\Delta t)e^{-r\Delta t}d\widehat{\mu }(x,y,t)\right\} .
\end{equation}%
Here $x$ belongs to grid $G$ and $y$ to $\mathcal{D}(x),\,$which\ is a set
of descendants of point $x,$ which may very well lie outside of the grid.
The function $\widehat{f}(y,t+\Delta t)$ is interpolated from the values of $%
\widehat{f}$ on the grid points that were obtained in the previous step of
recursion. Measure $d\widehat{\mu }$ approximates $d\mu $.

In summary, the main idea of the method is to separate two different uses of
the grid, which usually plays a crucial role in both approximating the
evolution of the stochastic process and in keeping information about the
option value function. We suggest using the grid only for the latter purpose
and simulating the stochastic process by computing small clusters of
descendant points around each grid point. The values on the descendant
points are interpolated from values on the grid points. Clearly, the success
of this idea crucially depends on the quality of the spatial interpolation.
We will discuss modern methods of multidimensional interpolation immediately
after presenting the outline of the algorithm and an example.

\section{Outline of Algorithm}

Here is the general outline of the algorithm:

1) Generate $G,$ a grid -- possibly irregular -- in the factor space.

2) For each point $x\in G,$ initialize the value function by computing the
payoff at the final stage $T$.

3) Begin recursion over $\,t<T$: Compute an approximation to the value
function at stage $t$ using an approximation at stage $t+1.$

This step can be realized in different ways. Since one part of this step is
applying the backward integral operator that corresponds to the factor
process, the algorithm becomes more precise if it can be computed
analytically. If not, we can always proceed as follows:

\qquad a) For each point $x\in G,$ compute a set of the states $\mathcal{D}%
(x)$ that follow $x$ in a discrete approximation to the factor process.

\qquad b) For each $y\in \mathcal{D}(x),$ compute $\widehat{f}(y,t+dt)$ by
interpolating.

\qquad c) Compute the continuation value function at point $x$ as the
discounted average of $\widehat{f}(y,t+dt)$ over $\mathcal{D}(x)$.

\qquad d) Compute the new value function by taking the maximum of the
continuation value function and the exercise payoff.

\qquad e) Proceed to the next step of the recursion.

The success of this method depends on the quality of approximations to the
continuation value fucntion. In multiple dimensions the most successful
approach to functional approximation to date is by fixing an over-complete
set of basis functions and then looking for an approximation recursively.
The typical realization of this idea is by the relaxed greedy algorithm
(RGA), which proceeds by forming a convex combination of a ridge function $%
\phi (b_{n+1}x+c_{n+1})$ and the previous approximation $f_{n}(x)$:%
\begin{equation}
f_{n+1}(x)=\alpha \phi (b_{n+1}x+c_{n+1})+(1-\alpha )f_{n}(x),
\end{equation}%
and then estimating $b_{n+1},c_{n+1},$ and $\alpha .$ Here $\phi (t)$ is a
function of one-dimensional parameter, and $b_{n+1}x$ means the scalar
product of vectors $b_{n+1}$ and $x.$

\citeN{barron93} advocates using a sigmoidal function $\phi $ consistent
with the neural network literature, \citeN{jones92} uses sinusoidal
functions, and \citeN{breiman93} suggests using a connected pair of
half-hyperplanes. The specifics of our application call for a different
choice. Since the approximation serves only as a tool for solving a partial
differential equation, the preferable choice seems to be Gaussian functions
that can be easily propagated backward by the action of the backward
integral operator associated with the factor process.

A question immediately arises: In which circumstances can recursive methods
find efficient approximations by Gaussians? Below is a theorem that answers
this question by generalizing a theorem from \citeN{jones92}.

Let 
\begin{equation}
f(x)=\sum_{i=1}^{\infty }a_{i}\phi (x;B_{i},C_{i})+u(x),
\label{function_type}
\end{equation}%
where $\phi (x;B,C)$ is a multi-dimensional Gaussian with parameters $B$
(precision, i.e. inverse of covariance matrix) and $C$ (shift of the center).

Assume that 
\begin{equation}
\begin{array}{cc}
1) & \sum \left| a_{i}\right| =K<\infty , \\ 
2) & \left\| \phi \right\| =1, \\ 
3) & \left\| u\right\| \leq \varepsilon ,%
\end{array}
\label{assumptions}
\end{equation}
where the norm is the $L^{2}$-norm.

Denote functions $\pm K\phi (x;B,C)$ that enter the expansion (\ref%
{function_type}) as $\phi _{i}$. Then $f$ is their linear convex
combination. Let us define the approximation at stage $n+1$ as a convex
combination $f_{n+1}=(1-\lambda )f_{n}+\lambda \phi $ of the previous stage
approximation and one of the functions $\phi _{i}$. We will choose the
approximation that solves the following problem 
\begin{equation}
\inf_{\lambda \in \lbrack 0,1],\phi \in \left\{ \phi _{i}\right\} }\left\{
\left\| (1-\lambda )f_{n}+\lambda \phi -f\right\| ^{2}\right\} .
\end{equation}

\begin{theorem}
\label{efficiency} For each $\alpha \in (0,1)$ and $n\geq 1$, the
approximant $f_{n}$ is guaranteed to satisfy the following inequality:%
\begin{equation}
\varepsilon _{n}^{2}\equiv \left\| f_{n}-f\right\| ^{2}\leq \max \left\{ 
\frac{\left[ (K+1)/\alpha \right] ^{2}}{n},\frac{\varepsilon ^{2}}{(1-\alpha
)^{2}}\right\} .
\end{equation}
\end{theorem}

Remark: By taking infinum over $\alpha $ we can obtain the following two
corollaries:

\begin{corollary}
If $\varepsilon =0,$ then $\varepsilon _{n}^{2}\leq (K+1)^{2}/n$
\end{corollary}

\begin{corollary}
If $\varepsilon >0$ and $n$ is sufficiently large$,$ then%
\begin{equation}
\varepsilon _{n}^{2}\leq \varepsilon ^{2}+2\frac{K+1}{\sqrt{n}}\varepsilon +%
\frac{(K+1)^{2}}{n}.
\end{equation}
\end{corollary}

\textbf{Proof of Theorem}: Let $f_{n}$ be the approximation obtained at step 
$n$ of the recursive algorithm. Consider the next recursive step. First, we
can write: 
\begin{eqnarray}
\left\| (1-\lambda )(f_{n}-f)+\lambda (\phi -f)\right\| ^{2} &=&(1-\lambda
)^{2}\left\| f_{n}-f\right\| ^{2}  \label{norm_combination} \\
&&+2\lambda (1-\lambda )\left( f_{n}-f,\phi -f\right) +\lambda ^{2}\left\|
\phi -f\right\| ^{2}.  \notag
\end{eqnarray}%
Next, using assumption (\ref{assumptions}.3) and the Cauchy-Schwarz
inequality, we obtain the following inequality:%
\begin{equation}
\left( f_{n}-f,\sum_{i=1}^{\infty }\alpha _{i}\phi _{i}-f\right) \leq
\left\| f_{n}-f\right\| \left\| \sum_{i=1}^{\infty }\alpha _{i}\phi
_{i}-f\right\| \leq \varepsilon _{n}\varepsilon ,
\end{equation}%
where $\alpha _{i}$ are coefficients of a convex linear combination.
Therefore, 
\begin{equation*}
\sum_{i=1}^{\infty }\alpha _{i}\left( f_{n}-f,\phi _{i}-f\right) \leq
\varepsilon _{n}\varepsilon
\end{equation*}%
From the positivity of $\alpha _{i}$ it follows that we can find such a $%
\phi _{i}$ that 
\begin{equation}
\left( f_{n}-f,\phi _{i}-f\right) \leq \varepsilon _{n}\varepsilon .
\end{equation}%
In addition, assumptions (\ref{assumptions}.1) and (\ref{assumptions}.2)
imply $\left\| \phi -f\right\| ^{2}\leq \left( \left\| \phi \right\|
+\left\| f\right\| \right) ^{2}\leq (K+1)^{2}$ for any $\phi .$
Consequently, (\ref{norm_combination}) implies 
\begin{eqnarray}
\varepsilon _{n+1}^{2} &\equiv &\inf_{\lambda \in \lbrack 0,1],\phi }\left\|
(1-\lambda )(f_{n}-f)+\lambda (\phi -f)\right\| ^{2} \\
&\leq &(1-\lambda )^{2}\varepsilon _{n}^{2}+2\lambda (1-\lambda )\varepsilon
\varepsilon _{n}+\lambda ^{2}(K+1)^{2}.
\end{eqnarray}%
Choose 
\begin{equation}
\lambda =\frac{(K+1)^{2}-\varepsilon \varepsilon _{n}}{(K+1)^{2}+\varepsilon
_{n}^{2}-2\varepsilon \varepsilon _{n}}.
\end{equation}%
This is a valid choice of $\lambda $ provided that $\varepsilon _{n}\geq
\varepsilon .$ For this choice we have the following bound for the error of
the next step:%
\begin{equation}
\varepsilon _{n+1}^{2}\leq \varepsilon _{n}^{2}\frac{(K+1)^{2}-\varepsilon
^{2}}{(K+1)^{2}+\varepsilon _{n}^{2}-2\varepsilon \varepsilon _{n}}.
\end{equation}%
Therefore, 
\begin{eqnarray}
\varepsilon _{n+1}^{-2} &\geq &\varepsilon _{n}^{-2}\left[ 1+\frac{\left(
\varepsilon _{n}-\varepsilon \right) ^{2}}{(K+1)^{2}-\varepsilon ^{2}}\right]
\\
&\geq &\varepsilon _{n}^{-2}+\frac{1}{(K+1)^{2}}\left( 1-\frac{\varepsilon }{%
\varepsilon _{n}}\right) ^{2}.
\end{eqnarray}%
Consequently, if $\varepsilon _{n}\geq \varepsilon /\left( 1-\alpha \right) $%
, then we have a sequence of inequalities:%
\begin{eqnarray}
\varepsilon _{n+1}^{-2}-\varepsilon _{n}^{-2} &\geq &\left( \frac{\alpha }{%
K+1}\right) ^{2}, \\
&&... \\
\varepsilon _{1}^{-2}-\varepsilon _{0}^{-2} &\geq &\left( \frac{\alpha }{K+1}%
\right) ^{2}.
\end{eqnarray}%
Summing them up we get 
\begin{eqnarray}
\varepsilon _{n+1}^{-2} &\geq &\varepsilon _{0}^{-2}+(n+1)\left( \frac{%
\alpha }{K+1}\right) ^{2}\geq (n+1)\left( \frac{\alpha }{K+1}\right) ^{2}, \\
\varepsilon _{n+1}^{2} &\leq &\frac{\left[ (K+1)/\alpha \right] ^{2}}{n+1}.
\end{eqnarray}%
Therefore, either $\varepsilon _{n}\leq \varepsilon /\left( 1-\alpha \right) 
$ and then $\varepsilon _{n+1}\leq \varepsilon /\left( 1-\alpha \right) $ $,$
or $\varepsilon _{n+1}^{2}\leq \left( n+1\right) ^{-1}\left( (K+1)/\alpha
\right) ^{2}.$

QED.

The significance of Theorem \ref{efficiency} is that it shows that a large
class of functions can be approximated to a high precision $\delta $ by an
expansion that has $O(\delta ^{-2})$ terms. In particular, the rate of
growth of the number of terms is independent of the space dimension.
Moreover, the theorem shows that the approximation can be found by the
recursive optimization method.

In practice the parameters of the expansion must be estimated from the
values of the function on a discrete grid. This introduces an additional
error in the approximation - the estimation error. The extent of this error
is not analysed in this paper. However, Theorem \ref{efficiency} and results
in \citeN{niyogi_girosi99} suggest that the number of the gridpoints needed
to bring the approximation error below a certain threshold grows only
polynomially with dimension.

\section{Lower and Upper Bounds}

In practical applications, we are often interested not only in an estimate
but also in the firm bounds on the option value. Fortunately, the Monte
Carlo method and interpolations from the lattice method can be combined for
the efficient calculation of these bounds.

Indeed, let $\widehat{V}(x,t)$ denote the approximation to the continuation
value function. Define the following stopping rule $\widehat{\tau }$: ``Stop
if $\widehat{V}(x,t)\leq \pi (x,t)$''. It simply tells the holder of the
option to stop when the approximate continuation value of the option is
smaller than the exercise value. Then the lower bound on the option value at
time 0 is given by 
\begin{equation}
V_{L}=E_{0}e^{-r\widehat{\tau }}\pi (x_{\widehat{\tau }},\widehat{\tau }%
)\leq \sup_{\tau }E_{0}e^{-r\tau }\pi (x_{\tau },\tau )=V.
\end{equation}%
The expectation can be easily computed using Monte Carlo simulations of
possible factor paths.

The upper bound can be computed using the duality method developed in %
\citeN{rogers02} and \citeN{haugh_kogan01}. Let $e^{-rt}M_{t}$ be a
supermartingale. Then 
\begin{eqnarray}
V &=&\sup_{\tau }E_{0}e^{-r\tau }\pi (x_{\tau },\tau ) \\
&=&\sup_{\tau }\left\{ E_{0}e^{-r\tau }\left[ \pi (x_{\tau },\tau )-M_{\tau }%
\right] +E_{0}e^{-r\tau }M_{\tau }\right\} \\
&\leq &\sup_{\tau }E_{0}e^{-r\tau }\left[ \pi (x_{\tau },\tau )-M_{\tau }%
\right] +M_{0}\text{ (because }E_{0}e^{-r\tau }M_{\tau }\leq M_{0}\text{)} \\
&\leq &E_{0}\sup_{t}e^{-rt}\left[ \pi (x_{t},t)-M_{t}\right] +M_{0}
\end{eqnarray}%
Consequently, if we manage to find such a supermartingale that the
expectation of $e^{-r\tau }\left[ \pi (x_{\tau },\tau )-M_{\tau }\right] $
is uniformly small, then we can calculate a good upper bound on the option
price.

Intuitively, supermartingale $M$ represents a replicating portfolio that the
writer of the option constructs to hedge his position. The portfolio should
be designed in such a way that it covers or almost covers the funds needed
in the case of the option exercise. The possible deficit in funds is
measured by the difference $\pi (x_{t},t)-M_{t}$, and the value of the
option cannot be larger than the sum of the replicating portfolio value $%
M_{0}$ and the discounted expectation of the supremum of the deficit.

The main question is how to find a good supermartingale $M.$ One possible
way is to choose the martingale distillation of the process $\overline{V}%
(x_{t},t)=e^{-rt}\max \left\{ \widehat{V}(x_{t},t),\pi (x_{t},t)\right\} .$
In this case we define 
\begin{equation}
e^{-r(t+1)}M_{t+1}=e^{-rt}M_{t}+E_{t}\overline{V}(x_{t+1},t+1)-\overline{V}%
(x_{t},t).
\end{equation}%
The expectation in this expression can be computed numerically. The choice
will provide a good upper bound provided that the approximation $\widehat{V}$
is good approximation to the option continuation value.

\section{Application}

The main goal of this section is to show that the method described above is
a viable alternative to approximate Monte Carlo methods. It demonstrates
this by computing values of rainbow options. These securities have the
following payoff structure:%
\begin{equation}
\pi =\max (f(K,S_{1},...,S_{n}),0).
\end{equation}%
For example, the put option on the minimum of several assets has 
\begin{equation}
f=K-\min (S_{1},...,S_{n}).
\end{equation}

This set of securities is a convenient benchmark for testing a pricing
method because it was extensively studied in the literature. For the case of
European puts on the minimum or the maximum there are analytic formulas
derived by \citeN{stulz82} and \citeN{johnson87}. The American put on the
minimum of two assets is priced in \citeN{boyle88} by a variant of the
multinomial lattice method. \citeN{boyle_evnine_gibbs89} give the results of
a lattice method for three-dimensional European puts. %
\citeN{broadie_glasserman97a}, \citeN{broadie_glasserman97b}, %
\citeN{raymar_zwecher97}, \citeN{broadie_glasserman_jain97}, and %
\citeN{longstaff_schwartz01} use two- and five-dimensional options as a
benchmark for their Monte Carlo simulation methods.

The most essential choice in our algorithm was how many gridpoints to use.
The grid was constructed by generating Sobol's quasi-Monte Carlo sequences
in a hypercube, and then transfoming them by an appropriate Gaussian
distribution. The number of grid points was determined using the
cross-validation method.\footnote{%
See \citeN{hastie_tibshirani01} for a fuller description of the
cross-validation method.} Namely, the grid was divided into two portions:
training and validation sets. The approximation was found using the training
set only, and its quality was evaluated on the validation set by computing
the mean squared error (MSE). The new terms in the approximative expansion
were added only if they decreased the MSE criterion. The lowest value of the
MSE was taken as the performance measure for a given number of data points.
If it was unsatisfactory, the number of grid points was increased and the
procedure repeated.

Table 1 summarizes the results of the application of the Interpolative
Lattice (IL) method to European and American puts on the minimum of two
assets. The results are compared to the results of the analytical formula
and to the results of the lattice method from \citeN{boyle88}. This table
shows that IL gives high precision results for European options where the
difference is on average less than a cent. For American puts the results of
the IL and Boyle methods are slightly different but still are very close to
each other.

Table 2 shows the results of IL application to valuation of various put
options on three assets. The parameters are as in %
\citeN{boyle_evnine_gibbs89} but since the results in %
\citeN{boyle_evnine_gibbs89} apparently contain a miscalculation, the
results obtained by the Monte Carlo method are used as a benchmark.\footnote{%
The Monte Carlo method has been tested on 2-dimensional options from %
\citeN{boyle88}, and the results have been found to be in remarkable
agreement with \ the results of \citeN{boyle88}.} The results for European
options are generally in good agreement with the Monte Carlo results, with
the difference less than $5\%$ of the option value.

Table 3 shows the results of pricing American puts on various functions of
five asset prices. The results are compared with results by the Least
Squares Monte Carlo, binomial and Berridge-Schumacher methods as they are
reported in \citeN{berridge_schumacher04}. The results are in good
agreement. For the put on the geometric average of the five asset prices,
the dimension reduction is possible and we are able to compute the exact
price. In this case , the interpolative lattice method slightly
overestimates the true price but the difference is less than 3\%. The
interpolative lattice method gives an estimate that in three cases exceeds
and in one case falls below the estimate given by Least Squares Monte Carlo.
This evidence suggests that the interpolative lattice method tends to give
an estimate that is biased upward.

The pricing times of the interpolative lattice method are rather long in the
case of five assets -- around 2 hours. Most of the time goes into the search
for a good adaptive approximation.

In summary, it appears that pricing by adaptive interpolation is a viable
alternative to other methods of pricing multidimensional options including
Least Squares Monte Carlo.

\section{\textbf{Conclusion}}

This paper proposes a novel variant of the lattice option pricing, which is
based on modern methods of adaptive interpolation. In multiple dimensions
the method allows using irregular grids and thus avoids both the curse of
dimensionality and the necessity to build recombining trees. The method is
easy to apply to many examples of derivative securities and its practical
viability is corroborated by numerical examples.

A theorem is demonstrated that suggests that the new method is less likely
to be vulnerable to the curse of dimensionality. Much further research is
needed, however, to determine the convergence properties of the new method.

\bibliographystyle{CHICAGO}
\bibliography{comtest}

\end{document}